\documentclass[12pt]{article}
\pdfoutput=1
\usepackage{natbib}
\usepackage{epsfig}
\usepackage{epsf}
\usepackage{here}
\usepackage{color}
\usepackage{pifont}
\newcommand{\DMdet}{\mbox{\boldmath ${\cal D}$}}

\newcommand{\dinf}[1]{\longrightarrow_{#1}}

\newcommand{\wli}[1]{\mbox{$\longrightarrow \hspace{-0.60cm} / \hspace{0.40cm}_{#1}~$}}

\newcommand{\FF}{{{\cal F}}}
\newcommand{\HH}{{{\cal H}}}

\newcommand{\X}{{{\bar X}}}

\newcommand{\hF}{{\bar F}}
\newcommand{\hG}{{\bar G}}
\newcommand{\hD}{{\bar D}}

\newcommand{\hf}{{\bar f}}

\newcommand{\IEf}{{\rm IEf}_{F\rightarrow D}}
\newcommand{\CEf}{{\rm CEf}_{F\rightarrow D}}
\newcommand{\PS}{{\cal S}}

\newcommand{\EI}{{\rm E}_{P^*_I}}
\newcommand{\EO}{{\rm E}_{P^*_O}}

\newcommand{\bX}{\mbox{\boldmath $X$}}

\newcommand{\bS}{\mbox{\boldmath ${\cal S}$}}

\newcommand{\CHD}{\mbox{\rm CHD}}
\newcommand{\BMI}{\mbox{\rm BMI}}
\newcommand{\betO}{\beta_{\mbox{\rm \tiny BMI}}}
\newcommand{\tBMI}{{\mbox{\rm \tiny BMI}}}
\newcommand{\tLDL}{{\mbox{\rm \tiny LDL}}}
\newcommand{\tCRP}{{\mbox{\rm \tiny CRP}}}
\newcommand{\tSBP}{{\mbox{\rm \tiny SBP}}}

\newcommand{\betL}{\beta_{\mbox{\rm \tiny LDL}}}
\newcommand{\betI}{\beta_{\mbox{\rm \tiny CRP}}}
\newcommand{\betP}{\beta_{\mbox{\rm \tiny SBP}}}

\newcommand{\LDL}{\mbox{\rm LDL}}
\newcommand{\CRP}{\mbox{\rm CRP}}
\newcommand{\SBP}{\mbox{\rm SBP}}

\newtheorem{Definition}{Definition}
\begin{document}

\title{The stochastic system approach to causality with a view toward lifecourse epidemiology}
\author{ Daniel Commenges$^{1,2}$ }
\maketitle
\noindent {\em $^{1}$INSERM, ISPED, Centre INSERM U-897-Epid\'emiologie-Biostatistique, Bordeaux, F-33000}\\
	{\em $^{2}$Univ. Bordeaux 2, ISPED, Centre INSERM U-897-Epid\'emiologie-Biostatistique, Bordeaux, F-33000}\\

\maketitle

{\bf Summary}. \vspace{3mm}

The approach of causality based on physical laws and systems proposed by \cite{Commenges2009} is revisited. The issue of "levels", the relevance to epidemiology and the definition of effects are particularly developed. Moreover it is argued that this approach that we call the stochastic system approach is particularly well fitted to study lifecourse epidemiology. A hierarchy of factors is described that could be modeled using a suitable multivariate stochastic process.  To illustrate this approach, a conceptual model for coronary heart disease mixing continuous and discrete state-space processes is proposed.


\section{Introduction}
There has been a growing interest in developing statistical formalisms which help to draw conclusions about causal influences from statistical analysis in various fields, and especially in epidemiology. Of course this raises the philosophical question about the very meaning of the concept of "cause", and here, different approaches have been taken. Thus different statistical formalisms arise from different philosophical views of causality. For general texts on causality we refer to \cite{bunge2008causality} and \cite{Pearl2000} for philosophical and statistical points of view respectively. In contemporary debates  three concepts have been proposed for grounding the concept of causality: intervention, counterfactuals and physical laws. Specific statistical methods have been developed on the basis of these concepts.

Intervention is an important concept because it may serve for defining causality. Intervention is the ultimate goal in many fields, and in particular in epidemiology. Randomized designs are considered as the best proof of a causal influence.
Defining causality on the basis of this concept however is not completely satisfactory. Moreover, randomized studies (and more generally experiments) are not sufficient and have limitations. The most obvious one is that sometimes experimentation is impossible or cannot be used to explore all aspects of the problem \citep{rosenblum2009analysing}.

 The counterfactual approach is an attempt to define causality. The concept has been defended by \cite{lewis1973causation} on the philosophical ground and first applied in statistical experiments by \cite{neyman1923application} and \cite{Rubin1978};  then it has been developed and applied in numerous papers in both experiments and observational epidemiology; for instance \cite{Robins2000} used this framework to develop the marginal structural models. Counterfactuals have been implemented in statistical models through the formalism of potential outcomes. This is nowadays the dominant approach in statistics. However defining causality through counterfactuals has raised many criticisms from a philosophical point of view \citep{sep-causation-counterfactual} and some in the statistical literature \citep{Dawid2000,Geneletti2010}. Moreover the potential outcome formalism seems difficult to apply in complex dynamical models.

The concept of physical laws, or laws of nature, may be used to ground causality from a philosophical point of view. The basic physical laws tell something about change in time, while time ordering is part of the definition of causality on which everybody agrees: causes precede effects. An interesting starting point is a famous text written by Laplace at the beginning of the XIX$^{th}$ century, which states that by knowing the laws of nature the future state of the universe could be predicted if the present state was known. Even if this program is somewhat too ambitious, it has been the aim of science to decipher the laws of nature, and prediction ability has been both a way for testing theories and and their main utility. More recently the concept of mechanism has been discussed by \cite{Wimsatt2007} and \cite{Bechtel2005}. In the statistical field this point of view has been translated in terms of dynamical models by \cite{Granger1969} using a time-series framework; \cite{Eichler2010} have developed general formula for computing marginal effects of interventions in time series framework. The approach has been generalized to continuous-time processes by \cite{Aalen1987}, \cite{Didelez2008} and \cite{Commenges2009}.


The aim of the paper is to focus the approach of \cite{Commenges2009} on the philosophical aspects and on the relevance in epidemiology, and especially lifecourse epidemiology, a field which attracts more and more interest \citep{Kuh2003}. It raises special challenges for causal interpretation, and implementation of the stochastic system approach appears to be well adapted to "understand trajectories", a main aim of lifecourse epidemiology \citep{Ben-Shlomo2007}.

The paper is organized as follows. In section 2 we give a philosophical account of the approach of causality as it can be grounded on physical laws applied to systems. In section 3 we present the general mathematical framework within which the application of physical laws to systems can be formalized: this is a general multivariate stochastic process, presented in \cite{Commenges2009}. Influences between components of the process can be defined through the Doob-Meyer decomposition and can be represented by a directed graph. The issue of quantifying an influence, that is the definition of effects is discussed. In section 4 we consider how this approach could be applied to lifecourse epidemiology, a topic of growing interest. We propose a classification of factors going from lifestyle factors to disease,  and we develop a conceptual model for coronary heart disease. Section 5 concludes.

\section{Causality}
\label{infl}
\subsection{Physical laws,  systems, levels}
Pierre-Simon Laplace wrote at the beginning of his Essai philosophique sur les probabilit\'es  \citep{laplace1825essai}:
"Nous devons donc envisager l'\'etat pr\'esent de l'univers, comme l'effet de son \'etat ant\'erieur, et comme la cause de celui qui va suivre. Une intelligence qui, pour un instant donn\'e, conna\^{\i}trait toutes les forces dont la nature est anim\'ee, et la situation respective des \^etres qui la composent, si d'ailleurs elle \'etait assez vaste pour soumettre ces donn\'ees \`a l'analyse, embrasserait dans la m\^eme formule les mouvements des plus grands corps de l'univers et ceux du plus l\'eger atome : rien ne serait incertain pour elle, et l'avenir comme le pass\'e serait pr\'esent \`a ses yeux. L'esprit humain offre, dans la perspective qu'il a su donner \`a l'astronomie, une faible esquisse de cette intelligence."
 In other (and english) words, if we knew the laws of nature (physical laws) and the state of the entities on which they apply at a given time, we could predict the future state of the universe. It is of course impossible to achieve such a program. The first obstacle is that we do not know all the physical laws and it is impossible to know the exact state of the universe at a given time. This was explicitly acknowledged by Laplace and was the reason why he went on developing an essay on probabilities just after formulating this plea for determinism ! Moreover quantum theory shows that randomness is an intrinsic feature of our universe and chaos theory shows how deterministic phenomena produce randomness.

 Another obstacle, less explicitly stated by Laplace, is that it is of course impossible to model the whole universe, not only because the precision of observation is limited but because we need manageable models. Only some very limited sets of events of the universe can be considered at a time. In the last sentence of the quoted text, Laplace takes as an illustration what has been achieved in predicting the movements of planets and comets using the recently developed theory of mass attraction. Such an analysis uses only Newton's laws applied to the movement of some celestial bodies identified only by their masses, and their locations and speeds at a given time. The fact that such an analysis is efficient shows that there is a property of "separability" of the universe. For instance, all the events relative to the presence of life on earth have no influence on the movement of earth; neither has the explosion of a star in a far-away galaxy. The development of a tumor in a human being is not affected by the movement of earth (conditionally on the fact that earth exists) and Newton laws are not very helpful to explain this phenomenon. This leads us to the concept of levels and of systems. Any scientific work begins by defining events of interest. These events are related to the evolution of entities (material or not). Physical laws allow predicting these events either exactly or, more generally, giving a probability law for  these events. The set of events and entities considered is called a system. An important concept, well known in the engineering literature \citep{sage1971system}, is that of "input", which makes the concept of system more flexible in that it is not completely isolated from the rest of the universe. Inputs can also serve for controlling the system.

  A system belongs to a certain {\em level}. Levels are indexed by scale and complexity. Scale is the dimension in space and time of the considered entities. For instance in astronomy we may consider the movement of planets of the solar systems during some centuries. Human beings do not belong to the spatial scale of this system. There is also a hierarchy of complexity. Elementary particles are structured in atoms, atoms in molecules, molecules in living cells, cells in organisms, organisms in societies. The scale and complexity dimensions interact in that considering a large scale is also considering a large number of units of a certain level. For instance considering a piece of metal of 10cm is considering a very marge number of atoms. Considering a region of space of one billion of cube-kilometers is considering such a large number of atoms that they can structure in stars and planets. Generally, to different levels, different physical laws apply. In physics, gravitation applies in principle at all scales but is negligible at the nanometer scale while it is the main force at the level of a planetary system. If we look at the complexity dimension, the gravitation law for instance is not very helpful for explaining the functioning of a living cell. The laws at a given level can in part be explained by the laws at a lower level. This reductionism principle is very useful, but at the same time there is the emergence of new laws at each level. This can be summarized by saying that "the whole is more than the sum of its parts" \citep{Anderson1972,Wimsatt1994}. We can distinguish systems made of a small number of entities (atoms, molecules, planetary systems) and those made of a large number of entities either homogeneous (crystal, piece of metal, gas, population of bacteria) or heterogenous (planet, cell, organism, society). In the first case we have a population of entities, in the second case we have a new entity at a different level. These dimensions of levels are illustrated in Figure 1 where some levels of different  complexity and scales are shown.\vspace{5mm}

\begin{figure}
\centering
\includegraphics[scale=0.6]{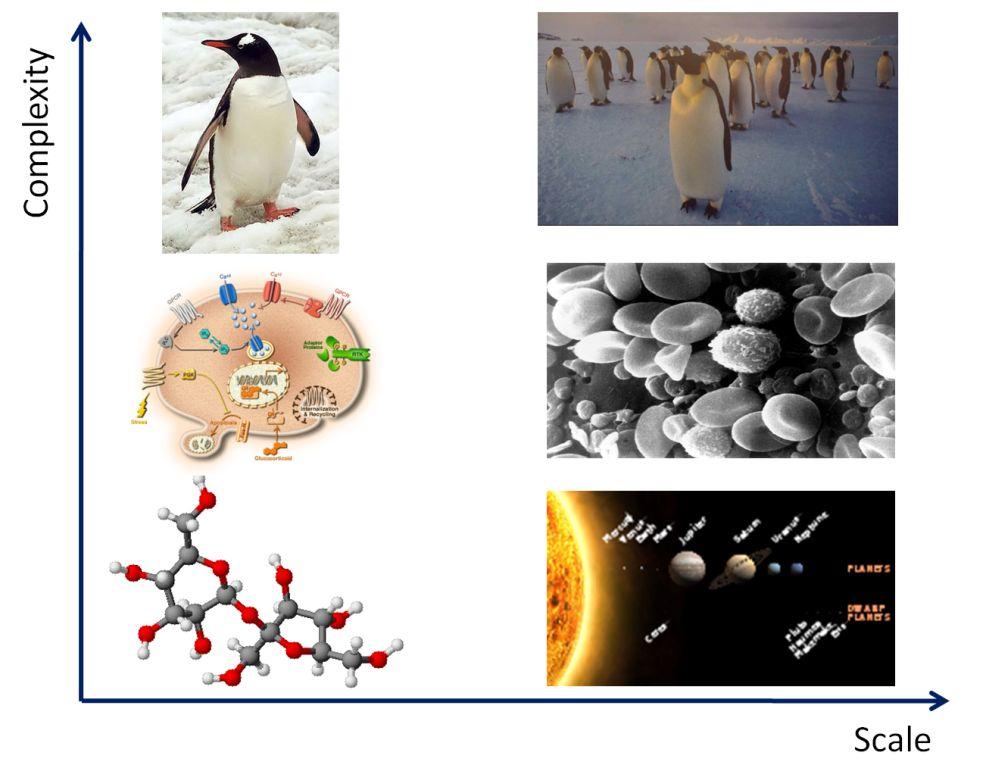}
\caption{Some levels. Vertical dimension: complexity: molecules, cells organisms are at different complexity levels. Horizontal dimension: scale: different levels arise if we consider large numbers of molecules, cells or organisms.}
\end{figure}

\subsection{Epidemiology as a science}

\label{causal-epid}
\subsection{Levels of Epidemiology}
In epidemiology we most often work at the level of a subject. Thus the systems are sets of events that can happen to a subject. As in the solar system example, we shall not consider all the events that are related to a subject, but rather selected sets of events that we consider as possibly connected, while being only weakly connected to any other event. For instance we may be interested in smoking and lung cancer. Note that we consider both smoking and lung cancer as "state" processes, and thus they will be represented by stochastic processes. On the contrary, gender, date of birth and genetic characteristics are attributes, that is not varying in time and contributing to the identification of that particular subject. As another example, we may consider HIV viral load, CD4+ T lymphocytes concentration and opportunistic diseases as components of the state, keeping gender, date of birth and genetic characteristics as attributes. While these factors are clearly attributes, this is less clear for other factors such as educational level, occupation, size. These factors could be considered as attributes in most studies because they are fairly constant on long time periods; they are however the result of processes which occurred during youth. Age is neither a process nor an attribute, but most often this is the time parameter in the processes associated to a subject.

While the central level of epidemiology is made of particular events happening to subjects (in general human beings) two other levels can be considered:  population of subjects (an upper level) and population of cells (a lower level). An epidemics is an event which happens at the population level. In pharmaco-epidemiology we may also be interested in the distribution of parameters in the population. Incidence and prevalence of a disease  are components of the "state" of a population. On the other hand, events occurring to populations of cells have an impact on clinical events at the subject level. For instance the interaction between populations of HIV viruses and CD4+ lymphocytes have an impact on AIDS. This level is at the margin of epidemiology and system biology.

\subsection{Laws in epidemiology}
\begin{figure}[h]

\centering
\includegraphics[scale=0.6]{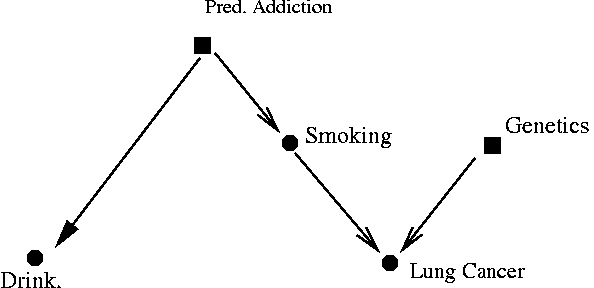}

\caption{Causal influences for Lung cancer. Filled circles represent states, while filled squares represent attributes. \label{lungcancer}}
\end{figure}

There is a huge gap between physics and epidemiology, with biology somewhat in-between. It may be questioned whether our approach of causality based on physical laws and inspired by physics is relevant here.  At a complex level we will speak rather of mechanisms than of physical laws (these mechanisms being based on physical laws at a lower level).
In epidemiology, in contrast with physics, very little is known about "physical laws". It may be questioned whether there are "physical laws" at the level of a human subject. Such laws could be called "clinical laws", reserving the term "epidemiological laws" to laws applying to populations of subjects. Such laws are more imprecise than in physics, they are most often probabilistic and they are not known with precision. Yet it is established that smoking increases the risk of lung cancer, and relatively precise studies have been done to describe the relationship between the two processes; for instance, \cite{breslow1987statistical}. In contrast there does not seem to be direct biological mechanisms which link alcohol drinking to lung cancer (although an indirect influence is not excluded). A graphical representation of causal influences for the factors ``Lung cancer'', ``Drinking'', ``Smoking'' could be constructed, adding a psychological predisposition to addictions that could be considered as an attribute as well as genetic factors (which are clearly attributes); see Figure \ref{lungcancer}.

Similarly, HIV infection has an influence on AIDS; the latency distribution between infection and AIDS has been estimated in several studies. These are clinical facts which derive from biological mechanisms. At the epidemiological level, such influences may be considered as direct, although at the biological level one may find intermediary factors. Concentrations of CD4+ T cells, of HIV viruses and of bacteria associated to an opportunistic disease are intermediary factors at the level of population of cells. These clinical facts can be used to study the effect of an epidemics in a population, for instance the relationship between incidence of infection and incidence of AIDS by backcalculation methods \cite{brookmeyer1991reconstruction}, or to predict the impact of an anti-smoking campaign on cancer incidence.
HIV infected patients treated with antiretroviral treatment have a much lower risk of developing AIDS than untreated patients,  opportunistic diseases may be prevented by prophylactic treatments. This is not really a law but a clinical fact that relies on biological mechanisms related to the properties of the treatment on the replication of the HIV virus and of the immune system. This fact lies on biological mechanisms, themselves lying on properties of molecules, themselves lying on basic physical laws.

\subsection{Learning laws in epidemiology}
At the beginning of any science both physical laws and relevant systems are unknown. The first step is the identification of entities. Then a precise description of events. Then a law is proposed and applied to a simple system. Then more complex systems are considered. The archetype of this process is astronomy of the solar system: identification of planets; description of their movement; Newton's laws; application to simple systems made of the Sun and one planet; more complex systems with several planets.
In epidemiology we are often in the situation of the beginning of a science where both laws and relevant systems are unknown. Medical doctors and epidemiologists however identify diseases, physiological conditions, behaviors and make progress in describing their evolutions. Approximate laws can sometimes be found which have a biological grounding. On the basis of these causal discoveries, interventions can be developed. Three big successes in which epidemiology has played an important role are the decrease of mortality from lung cancer, cardiovascular diseases and AIDS in most developed countries by intervening on both the main risk behaviors and the main biological factors. Note that intervening on behaviors implies considering an upper level of complexity that is psychology and sociology.

One obvious difficulty in epidemiology comes from the complexity of the biological processes leading to diseases. Another one is the sparsity of observations of these processes. Here the situation is improving as more and more cohort studies are launched with more and more repeated observations of biological markers (including imaging). One favorable feature in epidemiology is the availability of a large number of replica of similar systems (the events occurring in subjects of a large sample). Dynamical models are necessary to grasp the dynamics of the biological and clinical phenomena. Cox model was the first widely used approach attempting to grasp this dynamics. Systems are generally multivariate so that multistate models or joint models are more adapted to approach causality.

\section{Mathematical representations of laws and systems and definition of effects}

\subsection{Definition of influences}
\cite{Commenges2009} in line with \cite{Aalen1987} and \cite{Didelez2008} have developed a mathematical framework for representing causal influences. They defined a class of stochastic processes called $\DMdet$ which encompasses processes with continuous and discrete  state-spaces. For a given problem one defines a couple $(A,X)$, where $A$ is multivariate random variable representing the attributes of the system and $X$ is a multivariate stochastic process representing the state of the system. They further defined the WCLI property from a measurability condition of the compensator (or the intensity) of the process. Grossly speaking, a component $X_k$ is WCLI from a component $X_j$ if $X_j$ does not appear in the intensity of $X_k$. Influence is defined as a lack of the WCLI property. Consider for instance a two-dimensional stochastic process $X=(X_1,X_2)$ where $X_2$ admits the Doob-Meyer decomposition $dX_{2t}=\lambda_2(\X_{1t-},\X_{2t-})dt+dM_{2t}$, where $\X_{1t-},\X_{2t-}$ represent respectively the histories of $X_1$ and $X_2$ up to time $t$.  If $\lambda_2(\X_{1t-},\X_{2t-})$ does not depend on $\X_{1t-}$ (which can be given a rigorous definition in terms of measurability) then $X_2$ is WCLI from $X_1$. If $\lambda_2(\X_{1t-},\X_{2t-})$ does depend on $\X_{1t-}$, we shall say that $X_1$ directly influences $X_2$; we note $X_1 \dinf{\bX} X_2$. The definition of influence can be extended to influence of attributes on states.

 \cite{Commenges2009} also link the concepts of system, physical laws and causality; this is because the Doob-Meyer decomposition depends on both filtration and probability law. Essentially they assume that there exist a rich enough system $\bS^M$ in which influences between components of interests coincide with influences that could be computed using physical laws. In such a system influences are causal influences.

 WCLI can be defined on a random horizon $(0,T)$ suggested by \cite{Didelez2008} and \cite{Røysland2011}.
The extension to random horizon and intervals is useful in particular in applications to epidemiology. Assume that $X_k$ is a counting process representing the occurrence of a disease, lung cancer for instance. Let $T$ be the time of occurrence of lung cancer. We are interested in the influence of a risk factor, say smoking habits, on lung cancer. In our formalism lung cancer is represented by a counting process, and we wish to know whether the lung cancer process is WCLI from the risk factor process up to $T$. The influence of Lung cancer on the risk factor is another (generally less relevant) problem.

\subsection{Definition of conditional direct effects}
Most epidemiological studies  aim at assessing whether a factor $F$ has a causal influence on a disease $D$; if it does, it is important to quantify this influence. Here, we wish to go beyond the binary concept of causal influence: what we call ``effect'' is a degree of causal influence.

\begin{figure}[h]

\centering
\includegraphics[scale=0.6]{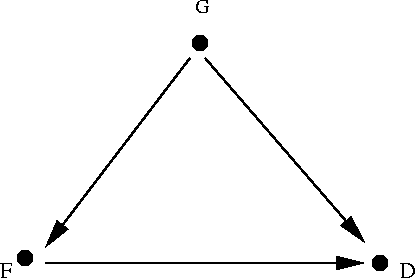}

\caption{"Perfect" observation system $\bS^M$.\label{simplesystem}}
\end{figure}

 We shall consider a simple situation, where we assume that there is a perfect system (for each subject in a population) that we call $\bS^M$, represented by three processes $F, D$ and $G$. We assume that $F$  causally influences $D$; $G$  causally influence both $F$ and $D$; there is no other causal influence on $D$. For simplicity we assume that there is no attribute. The graph of causal influences represented in Figure \ref{simplesystem} looks familiar, illustrating a confounding factor $G$, but there are some important differences with conventional graphs: (i) all the vertices represents processes; (ii) these could be multivariate processes; (iii) we assume that this is a ``perfect system'', in that the dynamics of the processes could in principle be computed from the knowledge of underlying biological mechanisms. The theory also allows reverse influence (the graph is not necessarily acyclic).

 The dynamics of $D$ can be analyzed using the Doob-Meyer decomposition (in the filtration $\FF^M$ generated by the process and under the true law $P^*$), which, assuming intensities exists, can be written as a system of stochastic differential equations (SDE):

 \begin{eqnarray} {dD_t} &=& \lambda_{Dt}(\hF_{t-},\hG_{t-},\hD_{t-})dt+ dM_{D,t} \label{eqD}\\
   {dG_t} &=& \lambda_{Gt}(\hG_{t-})dt+ dM_{G,t}\label{eqG}\\
{dF_t} &=& \lambda_{Ft}(\hF_{t-},\hG_{t-})dt+ dM_{F,t} \nonumber \end{eqnarray}

   It is natural to measure the direct conditional effect of $F$ as the contrast between  the values of the intensity of $D$  for two different values that $F$ could take; these values are in fact trajectories up to time $t-$, say $\hf_{t-}$ and $\hf'_{t-}$.
 We denote the instantaneous (causal) effect (or contrast) at time $t$, conditional on $G$ as:
\begin{equation} \IEf ^{G}(\hf_{t-},\hf'_{t-},t,\hD_{t-})= \phi(\lambda_{Dt}(\hf_{t-},\hG_{t-},\hD_{t-}),\lambda_{Dt}(\hf'_{t-},\hG_{t-},\hD_{t-})). \end{equation}
A cumulative effect can be defined in the of the cumulative intensities (or compensators) as :
\begin{equation} \CEf ^{G}(\hf_{t-},\hf'_{t-},t,\hD_{t-})= \phi(\Lambda_{Dt}(\hf_{t-},\hG_{t-},\hD_{t-}),\Lambda_{Dt}(\hf'_{t-},\hG_{t-},\hD_{t-})), \end{equation}
where $\phi(x,y)$ may be called a contrast function.

It is natural to take the contrast on either a multiplicative or an additive scale. Consider the case where $F$ and $G$ are univariate and not time-dependent, and $D$ is a disease represented by a counting process. If the intensity  $\lambda_{Dt}$ has the multiplicative form $\lambda_{Dt}=1_{\{D_{t-}=0\}}\alpha_{0}(t)e^{\beta_1 F+\beta_2 G}$, then the instantaneous effect is $e^{\beta_1 (f-f')}$, that is equal to the relative risk for a unit change of the value of $F$. In that particular case the instantaneous effect does not depend on time nor on the particular value of $G$; the cumulative effect takes the same value. On the other hand an additive contrast would fit better if the intensity itself has the additive form:
 \begin{equation} \label{additive} \lambda_{Dt}=1_{\{D_{t-}=0\}}[\alpha_{0}(t) + \beta_1 F+\beta_2 G].\end{equation}
 Then the instantaneous effect would be $\beta_1(f-f')$; the cumulative effect is $\beta_1(f-f')t$. In the additive model we have the equality: $\CEf ^G(f_{t-},f'_{t-},t)=\int_0^t \IEf ^G(f_{u-},f'_{u-},u)du$.

{\bf Remark.} [Effect depending on $D$] The effect as defined depends on the value of $D$ itself. Note that if $D$ is a 0-1 counting process, we can avoid mentioning the history of $D$, because up to the jump-time nothing happens. If $D$  has several jumps (representing recurrent events), for instance $D$ takes values $0,1, 2$, the effect of $F$ at time $t$ for $D_t=0$ and for $D_t=1$ may be different. This is not so surprising if we think that this system can be represented by a three-state irreversible process, and the effects of explanatory variables can be estimated for the two transitions $0\rightarrow 1$ and $1\rightarrow 2$.

In a general dynamical setting, there are advantages to assess the effects on an additive scale. Relative risks make sense if the intensities are positive; they make less sense if the ``intensities'' can be null or negative, which may happen in diffusion processes. Also, linearity of the contrast leads to simpler links between instantaneous and cumulative effects, and between conditional and mean effects.

Another issue is the definition of indirect effect. There is no space in this paper to develop this issue but there is an example in section \ref{indirecteffect}.

\subsection{Intervention systems}
In the previous section $F, G, D$ are generic names so that the system is abstract; in an application, we have a population of systems $\PS^{iM}$ corresponding to a population of subjects $i=1,\ldots,N$. The stochastic processes $(F^i,G^i,D^i)$ and $(F^{i'},G^{i'},D^{i'})$ are different processes, although they have the same clinical or biological meaning. For instance $D^i$ may represent lung cancer for subject $i$; it is not the same event if subject $i$ (say John) or subject $i'$ (say Paul) develops a lung cancer at time $t$. In some cases it is possible to construct systems in which $F$ can be manipulated, for instance for subjects in a clinical trial. The system for a subject $i'$ in a clinical trial will be specified by  $(F^{i'},G^{i'},D^{i'})$ with the same biological or clinical meaning as in an observation study; it will in addition contain a  control entry $U^{i'}$. This is also a perfect system for this subject but the true law governing $(F^{i'},G^{i'},D^{i'},U^{i'})$ is different from that of $(F^i,G^i,D^i)$. Since the way $D$ is influenced by $F$ derives from biological mechanisms, it is the same in both systems; the same is true for the influence of $G$ on $D$. However the presence of $U^{i'}$ makes the dynamics of $F$ different for $i$ and $i'$. If $F^{i'}$ is completely manipulated by $U^{i'}$ there can be no influence of $G^{i'}$ on $F^{i'}$.

{\bf Remark.} [Feasibility of intervention] One may ask how it is possible to construct systems in which $G^i\wli {\PS^{iM}} F^i$ leaving all the other influences unchanged. First it is not always possible to construct such systems. This is sometimes possible, and the reason why factor $F$ is of particular interest may be that $F$ can be easily manipulated. We may consider that in such cases the influence of $G$ on $F$ is ``fragile'' and can be broken; this is the case if $F$ is a treatment.

\begin{figure}[h]

\centering
\includegraphics[scale=0.6]{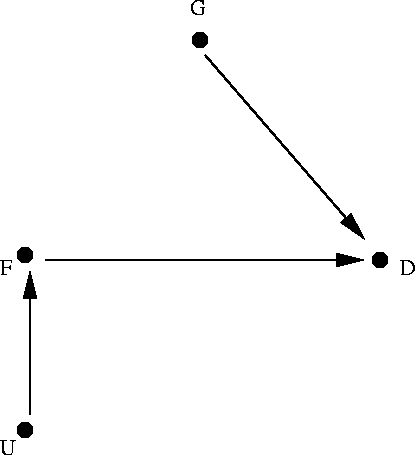}

\caption{"Perfect" intervention system $\PS^{Ii'M}$.\label{intervention}}
\end{figure}

Mathematically we could say that an observational system $\PS^{iM}$ and an interventional system $\PS^{Ii'M}$ are the same but with a different law. These laws are the restriction of the true probability law $P^*$ to the filtrations of the systems $\PS^{iM}$ and $\PS^{Ii'M}$, that we will call $P^*_O$ and $P^*_I$ respectively, omitting the subscripts $i$ and $i'$ in the iid case. However we must not forget that these systems are different as they correspond to different subjects. A generic system, where the stochastic processes do not represent real events, but  just retains the structure of possible real systems will be called an abstract system. In the following we shall continue with abstract systems for notational simplicity. We can now study the relationship between an abstract observation system (Figure \ref{simplesystem}) and a corresponding abstract intervention system (Figure \ref{intervention}). One goes from the first to the second by breaking the influence of $G$ on $F$ while keeping the other influences unchanged, an operation that \cite{Pearl2000} represents by the ``Do'' operator. It is clear that the conditional effect of $F$ on $D$ is the same in both models because $\lambda_{Dt}(\hF_{t-},\hG_{t-},\hD_{t-})$ is the same.

\subsection{Computation of marginal causal effect in intervention and observation systems}
\subsubsection{General}
An important issue is to be able to compute the mean effect of an intervention (the marginal causal effect).
By the innovation theorem, the intensity of $D$ in the simpler intervention system $\PS^{I2}$ where $G$  is ignored, is:  $\EI [\lambda_{Dt}(\hF_{t-},\hG_{t-},\hD_{t-})|\hF_{t-},\hD_{t-}]$.
This is interesting because we may not always observe $G$ and we may wish to have a summary of the effect of $F$ not depending on $G$. This is often not very easy to compute mathematically but this will be directly observed in the intervention system $\PS^{I2}$.

\begin{figure}[h]
\centering
\includegraphics[scale=0.7]{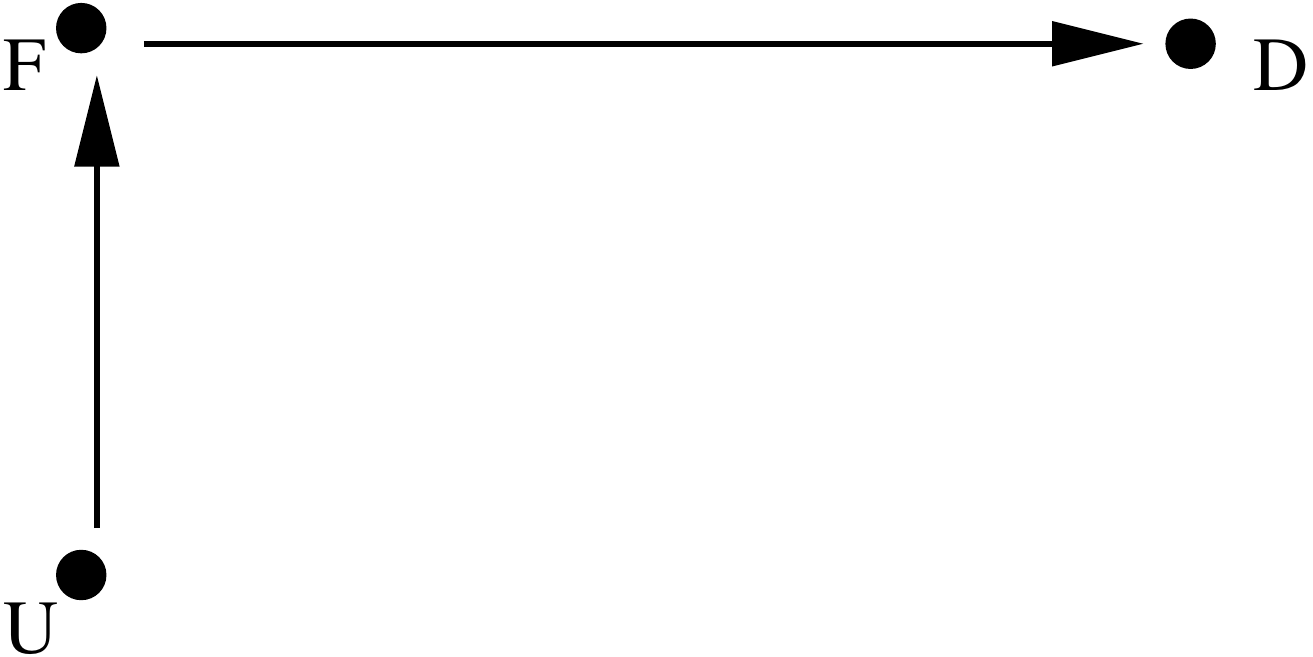}
\caption{Imperfect intervention system $\PS^{I2}$\label{impintsys}}
\end{figure}

In contrast the marginal causal effect cannot in general be obtained if we know only the law of the simpler observational system $\PS^{2}$ which does not contain $G$ (intuitively, this means that we cannot estimate the effect of $F$ if we don't take confounders into account). In such a system, by the innovation theorem, the compensator of $D$ in the filtration $(\FF^2)$ generated by $F$ and $D$ is $\lambda_{Dt}^{\FF^2}(\hF_{t-})={\EO}[\lambda_{Dt}(\hF_{t-},\hG_{t-},\hD_{t-})|\FF^2_{t-}]$. This expectation is not equal in general to ${\EI} [\lambda_{Dt}(\hF_{t-},\hG_{t-},\hD_{t-})|\FF^2_{t-}]$.

If we observe the perfect observation system it is possible to retrieve the marginal causal effect. This is because if we know $P_O^*$ on $\PS^M$ we also know $P_I^*$ (see section \ref{simplexample}). In practice if $G$ and $D$ are complex processes this may not be so easy. It would still be possible to compute any relevant effect by simulation. Indeed the system $\PS^M$ is characterized under $P_I^*$ where $\hF=\hf$ by the Doob-Meyer decomposition:

\begin{eqnarray*} {dD_t} &=& \lambda_{Dt}(\hf_{t-},\hG_{t-},\hD_{t-})dt+ dM_{D,t}\\
    {dG_t} &=& \lambda_{Gt}(\hG_{t-},)dt+ dM_{G,t}.\end{eqnarray*}

These equations are the same as the first two equations (\ref{eqD}),(\ref{eqG}) for $P_O^*$. Realizations of this system of SDE can be generated by simulation:  if both $G$ and $D$ are counting processes this is not difficult, if they are diffusion processes this can be done using techniques described in \cite{kloeden1994stochastic}; these techniques must be adapted if there is a mix of counting and diffusion processes. Realizations of $D$ allow to estimate its marginal law.


\subsubsection{A simple example}\label{simplexample}
Let us consider the case where $D$ is a disease represented by a counting process. If the processes $F$ and $G$ are time-constant they can be summarized by random variables. It is then equivalent to specify the law of $D$ by the intensity of the process or by the survival distribution $S(t|f)=P(D_t=0|F=f)$, and we have:  $S(t|f)=P(G=1|F=f) S(t|f,G=1)+P(G=0|F=f) S(t|f,G=0)$. In the intervention system, that is under the law $P_I^*$, $G$ and $F$ are marginally independent. Thus we have $P_I^*(G=1|F=f)=P_I^*(G=1)$. It follows that
\begin{equation} \label{MCE} S_I(t|f)=P(G=1) S(t|f,G=1)+P(G=0) S(t|f,G=0)\end{equation}
 This is what is observed (using a large sample) in the imperfect intervention system  represented in Figure \ref{impintsys}.
$S(t|f,G=1)$ and $S(t|f,G=0)$ (the survival conditional on both $F$ and $G$) as well as $P(G=1)$ and $P(G=1)$ (the marginal distribution of $G$) are the same in the observation and intervention systems.
Thus if we know the dynamics of the perfect observation system we can compute the marginal causal effect given by formula (\ref{MCE}). Note that we can also compute the marginal intensity of $D$; here $G$ has the role of a random effect and this problem has been studied by \cite{aalen2008survival} (section 6.2).

In contrast when we use the imperfect observation system depicted in Figure \ref{impobssys} what we observe is $$S_O(t|f)=P(G=1|F=f) S(t|f,G=1)+P(G=0|F=f) S(t|f,G=0),$$ and we cannot reconstitute $S_I(t|f)$ because we observe neither $P(G=1)$ nor $S(t|f,G=g)$ for $g=0,1$.

As an applied example in epidemiology, $F$ could represent ``Drinking'' and $G$  ``Smoking'': the link between smoking and drinking could produce a spurious influence of drinking on lung cancer if we do not adjust on smoking.

\section{Conceptual models for lifecourse epidemiology}
\subsection{The stochastic system approach to lifecourse epidemiology}
The stochastic system approach could be applied to lifecourse epidemiology. In practice this entails defining processes which represent how physiological processes evolve during life. The main difficulty is to define a relatively small number of processes that are sufficiently meaningful to make a good system. Another major difficulty arises when we consider the whole life from birth (or even from conception): some features evolve fast during childhood but stabilize during adulthood. The most obvious one is height, which can be considered as an attribute for adults (this is given on identity cards) but would be considered as part of the state for a child. So, we will restrict here to adult life epidemiology. It is useful to distinguish five categories of processes:
 \begin{enumerate}
\item {\bf lifestyle factors}: diet, physical activity, drinking, smoking, sexual behavior,...,
\item {\bf physiological conditions} (or "embodied" processes): body mass, lipid profile, blood pressure, inflammation level, oxidative stress, HIV infection,...,
\item {\bf pathological processes}: atheromatous process, pre-cancerous cells, neuronal degeneration, loss of CD4+ T-lymphocytes,...,
\item {\bf diseases}: coronary heart disease, cancers, Alzheimer disease, AIDS,...,
\item {\bf death}.
\end{enumerate}

In order to construct a model one has first to choose indicators of these broadly defined processes, especially for physiological conditions; for instance body mass can be represented by body mass index (BMI), lipid profile by LDL cholesterol concentration, high blood pressure by systolic blood pressure (SBP), and so on. Pathological processes are often difficult to observe and may be represented by latent processes. Diseases are generally well classified and observed. Death is of course the simplest process to define and observe. As for the lifestyle factors, one must also define indicators and it is often difficult to observe them accurately. Generally as we progress in the list, the conditions are less and less reversible, and death is not reversible at all. On the other hand lifestyle factors are themselves influenced by psycho-sociological factors. In order to avoid the system to grow indefinitely we may consider lifestyle factors as input to the system. Finally genetic factors, which are attributes of the subjects, may influence all the processes, although with a more modest impact than is generally thought.

\subsection{A conceptual model for adult life epidemiology of coronary heart disease}
\subsubsection{The model for coronary heart disease}
\begin{figure}[h]

\centering
\includegraphics[scale=0.6,angle=-90] {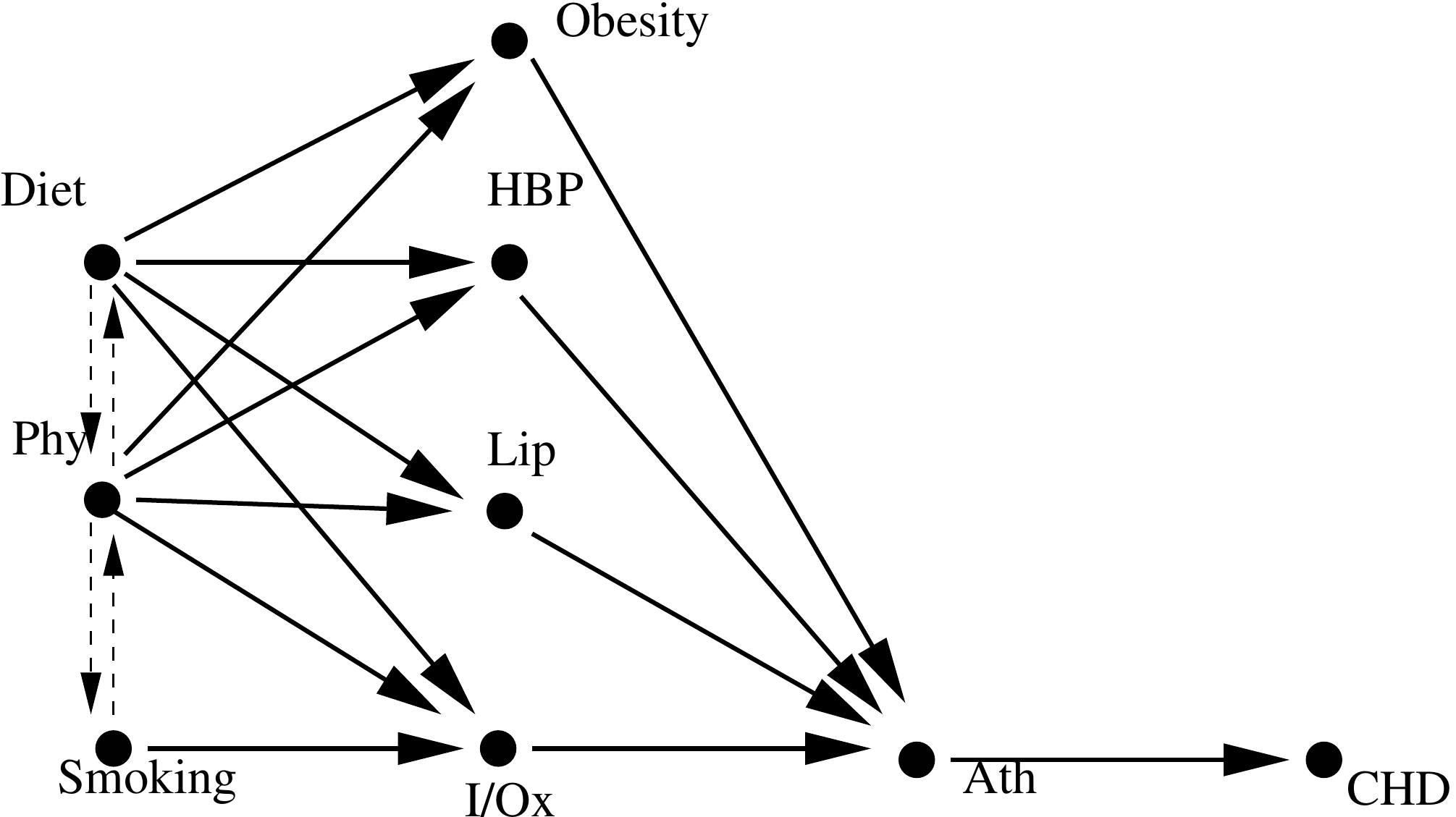}

\caption{Causal pathway for smoking and CHD: Phy: Physical activity; HBP: high blood pressure; Lip: lipids profile (cholesterol, low and high density,...); I/Ox: process representing a complex of inflammatory process and oxidative stress; Ath: atheromatous process.\label{Smoking-CHD}}
\end{figure}

We focus the above approach to coronary heart disease (CHD). Atherosclerosis is the main mechanism of CHD \citep{Nicholls2006} and the main risk factors of CHD are known: \cite{Berliner1995,Castelli1996,Ridker1998,Kannel2002,Emberson2003,VanGaal2006}.
Beyond a mere enumeration of risk factors our approach consists in constructing a causal pathway between processes.
The "lifestyle" factors are diet, physical activity, smoking; the physiological conditions are obesity, lipid profile, oxydative stress, inflammatory processes; the pathological process is the atheromatous process and the disease is CHD, which includes myocardial infarction and may lead to death. On this physiological basis we may attempt to construct a model, as simple as possible, but reflecting the dynamics of these processes and the pathways leading from lifestyle to disease. A possible conceptual model is shown in Figure \ref{Smoking-CHD}.

To go further we have to represent the categories like "obesity", "blood pressure", and so on, by quantitative indicators, and the set of values taken by them will form a stochastic process. Choosing binary indicators would lead to a multistate model with many states, that can also be represented by a multivariate counting process such as proposed by \cite{commenges2007likelihood}. For instance a good indicator for obesity is the BMI; it is possible to define obesity as $\BMI\ge 30$. Similarly hypercholesterolemia, hypertension and so on can be defined by dichotomization. This however entails a loss of information. Most of the processes are better represented as having a continuous state space. In order to obtain manageable systems of equations we assume that, possibly after suitable transformations, conditionally on the lifestyle factors, the physiological processes may be modeled by processes which have the same structure as an Ornstein-Uhlenbeck (OU) processes; see \cite{kallenberg2002foundations} for theory and \cite{zhu2011stochastic} for recent use in  biostatistics.
For instance  we assume that (maybe after transformation) the LDL process can be described by the following stochastic differential equation:

\begin{equation} \label{SDE-OU}
d\LDL_t=\theta_{\tLDL}(\LDL_t-\mu_{\tLDL,t})dt+dB_{\tLDL ,t},\end{equation}
where $B_{\tLDL ,t}$ is a Brownian motion. $\mu_{\tLDL,t}$ is modelled as a function of indicators of the lifestyle factors "total energy intake" for Diet ($D_t$), "total energy expenditure" for Physical activity ($\phi_t$), "number of cigarettes per day" for smoking ($S_t$) which are considered as input functions. The simplest model for $\mu_{\tLDL,t}$ is the linear model:
$$ \mu_{\tLDL,t}=\mu_{0,\tLDL}+ \beta_{S,\tLDL}S_t+\beta_{\phi, \tLDL}\phi_t+\beta_{D,\tLDL}D_t.$$
Note that if the lifestyle factors change in time, $ \mu_{\tLDL,t}$ is a moving target; if they are constant, $\LDL$ is an ordinary OU process and gets close to a stationary process.
For indicators of obesity, high blood pressure and inflammation $\BMI_t,\SBP_t,\CRP_t$ similar modeling can be proposed.
A diffusion process with a positive intensity is adapted to reflect the accumulative character of the atheromatous process.
\begin{equation} dA_t=\lambda_{At}dt+dB_{At}, \end{equation}
where the intensity can be modeled as a linear function of the physiological conditions:
\begin{equation} \lambda_{At}=\lambda_0+\betO \BMI_t+\betL \LDL_t+\betI \CRP_t +\betP \SBP_t. \end{equation}

Finally CHD can be modeled as an event which occurs when the atheromatous process crosses a certain threshold: $\CHD_t=1_{\{A_t>\eta\}}$. Thus, for CHD we assume a degradation model as used by \cite{Whitmore1986} and \cite{aalen2008survival} among others.

\subsubsection{Defining effects}\label{indirecteffect}
Putting aside the inference issue, that is pretending that there are values of $\beta$ and $\eta$ which determine the true probability under which the events are produced and that there are known, the question that is addressed is how to quantify the effect of a lifestyle process, for instance "smoking" on the atheromatous process and on myocardial infarction. This is an indirect effect, mediated by the physiological conditions. The first step is to express the process of interest, for instance $A_t$ as a function of the lifestyle processes. For constant lifestyle processes this can be done analytically because the solution of the differential equation (\ref{SDE-OU}) is an OU process.
For instance the LDL process is:
\begin{equation} \label{eqLDL} \LDL_t=\LDL_0 e^{\theta_{\tLDL} t}+\mu_{\tLDL,t}(1-e^{-\theta_{\tLDL} t})+\int_0^te^{\theta_{\tLDL} (s-t)}dB_{\tLDL,s}.\end{equation}
Taking into account the linear model of $\mu_{\tLDL,t}$ and assuming the initial values of the embodied processes are equal to zero, the intensity of the atheromatous process can be written:
\begin{eqnarray*}\lambda_{At}&=&\lambda_0+[\betO^t \beta_{S,\tBMI}+\betL^t \beta_{S,\tLDL}+\betI^t \beta_{S,\tCRP}+\betP^t \beta_{S,\tSBP}]S_t\\
        &+&[\betO^t \beta_{\phi,\tBMI}+\betL^t \beta_{\phi,\tLDL}+\betI^t \beta_{\phi,\tCRP}+\betP^t \beta_{\phi,\tSBP}]\phi_t\\
        &+&[\betO^t \beta_{D,\tBMI}+\betL^t \beta_{D,\tLDL}+\betI^t \beta_{D,\tCRP}+\betP^t \beta_{D,\tSBP}]D_t +M_t,
         \end{eqnarray*}
where $M_t$ is a martingale coming from the last term of equation (\ref{eqLDL}) and other similar equations,  $\betL^t=\betL (1-e^{-\theta_{\tLDL} t})$ and similar definitions for $\betL^t,\betI^t,\betP^t$.
Up to a term of zero expectation, a linear instantaneous  contrast between two trajectories of the function $S$, $s$ and $s'$ for the total effect on the atheromatous process gives:
\begin{equation}\lambda_{At}(s)-\lambda_{At}(s')=[\betO^t \beta_{S,\tBMI}+\betL^t \beta_{S,\tLDL}+\betI^t \beta_{S,\tCRP}+\betP^t \beta_{S,\tSBP}](s_t-s'_t).\end{equation}
Such a simple formula however is not available for summarizing the effect of smoking on CHD. Customized contrasts can be computed, for instance by computing $P(\CHD_t=1|{\rm Do }~ S_t=s_t)$.

\subsubsection{Possible implementation}
If the lifestyle processes are controlled, the atheromatous process is a Gaussian process. The distribution of the time for hitting the barrier $\eta$ could be computed numerically by simulation. Thus the computation of the likelihood for such a model is challenging. There is also the problem of finding enough information for estimating the parameters. It is likely that one cannot find a single study giving enough information to identify the model and that we have to resort n some kind of synthesis analysis \citep{presanis2011bayesian}.

With the assumption that the lifestyle processes are constant, that the equilibrium is reached rapidly and neglecting the variability, the physiological processes become constant and the atheromatous process is a Brownian process with linear drift with slope $\lambda=\lambda_{At}$. In that case the distribution of the hitting time is an Inverse Gaussian distribution. Its  p.d.f. is $\frac{\eta}{\sqrt{2\pi t^3}}e^{-\frac{(\eta-\lambda t)^2}{2t}}$, so that all the probabilities of all the events pertaining to CHD can be rather easily computed. There is however no simple form for a contrast between different values of a lifestyle process. For given values of the parameters and values of the lifestyle processes the hazard and the survival functions can be computed and visualized.
Such a model has been fitted to the data of the literature by \cite{Commenges2012Evidence}. They used a synthesis approach to fit heterogenous results of several studies by maximizing a pseudo-likelihood. They obtained reasonably good fit and consistent values of regression coefficients.

\section{Conclusion}

We have given an account of a philosophical approach to causality based on physical laws and systems and we have recalled the main mathematical tools used for representing causal relationships and proposed some ways of computing effects. We have examined how this could be applied to epidemiology, and particularly to lifecourse epidemiology. In practice this leads to develop dynamical models based on stochastic processes, incorporating as much science as possible. There are two obstacles for implementing complex dynamical models for lifecourse epidemiology. In complex models the distribution of the observable cannot be obtained analytically and thus they must be approximated by simulations. This seems however feasible, especially with the development of parallel computing. The other main obstacle is to gather sufficient information to estimate the parameters for all the stages of the model. There are however more and more large or very large cohorts with a long follow-up. Moreover the synthesis approach used by \cite{Commenges2012Evidence}, see also \cite{presanis2011bayesian}, is a promising way for identifying the models by finding information about diverse parameters in a large number of studies.

\section{Acknowledgements}
I thank Anne G\'egout-Petit for stimulating discussion about the topic of this paper.

\bibliographystyle{chicago}
\bibliography{Causal}

\begin{thebibliography}{}

\bibitem[\protect\citeauthoryear{Aalen}{Aalen}{1987}]{Aalen1987}
Aalen, O. (1987).
\newblock {Dynamic modelling and causality}.
\newblock {\em Scandinavian Actuarial Journal\/}.

\bibitem[\protect\citeauthoryear{Aalen, Borgan, and Gjessing}{Aalen
  et~al.}{2008}]{aalen2008survival}
Aalen, O., {\O}.~Borgan, and H.~Gjessing (2008).
\newblock {\em Survival and event history analysis: a process point of view}.
\newblock Springer Verlag.

\bibitem[\protect\citeauthoryear{Anderson}{Anderson}{1972}]{Anderson1972}
Anderson, P.~W. (1972, August).
\newblock {More is different.}
\newblock {\em Science (New York, N.Y.)\/}~{\em 177\/}(4047), 393--6.

\bibitem[\protect\citeauthoryear{Bechtel and Abrahamsen}{Bechtel and
  Abrahamsen}{2005}]{Bechtel2005}
Bechtel, W. and A.~Abrahamsen (2005, June).
\newblock {Explanation: a mechanist alternative}.
\newblock {\em Studies in History and Philosophy of Science Part C: Studies in
  History and Philosophy of Biological and Biomedical Sciences\/}~{\em
  36\/}(2), 421--441.

\bibitem[\protect\citeauthoryear{Ben-Shlomo}{Ben-Shlomo}{2007}]{Ben-Shlomo2007}
Ben-Shlomo, Y. (2007, June).
\newblock {Rising to the challenges and opportunities of life course
  epidemiology.}
\newblock {\em International journal of epidemiology\/}~{\em 36\/}(3), 481--3.

\bibitem[\protect\citeauthoryear{Berliner, Navab, Fogelman, Frank, Demer,
  Edwards, Watson, and Lusis}{Berliner et~al.}{1995}]{Berliner1995}
Berliner, J.~A., M.~Navab, A.~M. Fogelman, J.~S. Frank, L.~L. Demer, P.~A.
  Edwards, A.~D. Watson, and A.~J. Lusis (1995, May).
\newblock {Atherosclerosis: basic mechanisms. Oxidation, inflammation, and
  genetics.}
\newblock {\em Circulation\/}~{\em 91\/}(9), 2488--2496.

\bibitem[\protect\citeauthoryear{Breslow, Day, et~al.}{Breslow
  et~al.}{1987}]{breslow1987statistical}
Breslow, N., N.~Day, et~al. (1987).
\newblock Statistical methods in cancer research. volume ii--the design and
  analysis of cohort studies.
\newblock {\em IARC scientific publications\/}~(82), 1.

\bibitem[\protect\citeauthoryear{Brookmeyer}{Brookmeyer}{1991}]{brookmeyer1991%
reconstruction}
Brookmeyer, R. (1991).
\newblock Reconstruction and future trends of the aids epidemic in the united
  states.
\newblock {\em Science\/}~{\em 253\/}(5015), 37.

\bibitem[\protect\citeauthoryear{Bunge}{Bunge}{2008}]{bunge2008causality}
Bunge, M. (2008).
\newblock {\em Causality and modern science}.
\newblock Transaction Pub.

\bibitem[\protect\citeauthoryear{Castelli}{Castelli}{1996}]{Castelli1996}
Castelli, W.~P. (1996, July).
\newblock {Lipids, risk factors and ischaemic heart disease}.
\newblock {\em Atherosclerosis\/}~{\em 124\/}(1), S1--S9.

\bibitem[\protect\citeauthoryear{Commenges and G{\'e}gout-Petit}{Commenges and
  G{\'e}gout-Petit}{2007}]{commenges2007likelihood}
Commenges, D. and A.~G{\'e}gout-Petit (2007).
\newblock Likelihood for generally coarsened observations from multistate or
  counting process models.
\newblock {\em Scandinavian journal of statistics\/}~{\em 34\/}(2), 432--450.

\bibitem[\protect\citeauthoryear{Commenges and G\'{e}gout-Petit}{Commenges and
  G\'{e}gout-Petit}{2009}]{Commenges2009}
Commenges, D. and A.~G\'{e}gout-Petit (2009, June).
\newblock {A general dynamical statistical model with causal interpretation}.
\newblock {\em Journal of the Royal Statistical Society: Series B (Statistical
  Methodology)\/}~{\em 71\/}(3), 719--736.

\bibitem[\protect\citeauthoryear{Commenges and Hejblum}{Commenges and
  Hejblum}{2012}]{Commenges2012Evidence}
Commenges, D. and B.~Hejblum (2012).
\newblock Evidence synthesis through a degradation model applied to myocardial
  infarction.
\newblock {\em submitted\/}.

\bibitem[\protect\citeauthoryear{Dawid}{Dawid}{2000}]{Dawid2000}
Dawid, A. (2000).
\newblock {Causal inference without counterfactuals}.
\newblock {\em Journal of the American Statistical Association\/}~(188),
  407--424.

\bibitem[\protect\citeauthoryear{Didelez}{Didelez}{2008}]{Didelez2008}
Didelez, V. (2008, January).
\newblock {Graphical models for marked point processes based on local
  independence}.
\newblock {\em Journal of the Royal Statistical Society: Series B (Statistical
  Methodology)\/}~{\em 70\/}(1), 245--264.

\bibitem[\protect\citeauthoryear{Eichler and Didelez}{Eichler and
  Didelez}{2010}]{Eichler2010}
Eichler, M. and V.~Didelez (2010).
\newblock {On Granger causality and the effect of interventions in time
  series}.
\newblock {\em Lifetime data analysis\/}~{\em 16}, 3--32.

\bibitem[\protect\citeauthoryear{Emberson}{Emberson}{2003}]{Emberson2003}
Emberson, J. (2003, October).
\newblock {Re-assessing the contribution of serum total cholesterol, blood
  pressure and cigarette smoking to the aetiology of coronary heart disease:
  impact of regression dilution bias}.
\newblock {\em European Heart Journal\/}~{\em 24\/}(19), 1719--1726.

\bibitem[\protect\citeauthoryear{Geneletti and Dawid}{Geneletti and
  Dawid}{2010}]{Geneletti2010}
Geneletti, S. and A.~Dawid (2010).
\newblock {Defining and identifying the effect of treatment on the treated}.
\newblock {\em Research Report\/}~{\em 3}, 1--24.

\bibitem[\protect\citeauthoryear{Granger}{Granger}{1969}]{Granger1969}
Granger, C. W.~J. (1969, August).
\newblock {Investigating Causal Relations by Econometric Models and
  Cross-spectral Methods}.
\newblock {\em Econometrica\/}~{\em 37\/}(3), 424.

\bibitem[\protect\citeauthoryear{Kallenberg}{Kallenberg}{2002}]{kallenberg2002%
foundations}
Kallenberg, O. (2002).
\newblock {\em Foundations of modern probability}.
\newblock Springer Verlag.

\bibitem[\protect\citeauthoryear{Kannel, Wilson, Nam, and D'Agostino}{Kannel
  et~al.}{2002}]{Kannel2002}
Kannel, W.~B., P.~W.~F. Wilson, B.-H. Nam, and R.~B. D'Agostino (2002,
  October).
\newblock {Risk stratification of obesity as a coronary risk factor.}
\newblock {\em The American journal of cardiology\/}~{\em 90\/}(7), 697--701.

\bibitem[\protect\citeauthoryear{Kloeden, Platen, and Schurz}{Kloeden
  et~al.}{1994}]{kloeden1994stochastic}
Kloeden, P., E.~Platen, and H.~Schurz (1994).
\newblock {\em Numerical Solution of SDE Through Computer Experiments}.
\newblock Springer.

\bibitem[\protect\citeauthoryear{Kuh}{Kuh}{2003}]{Kuh2003}
Kuh, D. (2003, October).
\newblock {Life course epidemiology}.
\newblock {\em Journal of Epidemiology \& Community Health\/}~{\em 57\/}(10),
  778--783.

\bibitem[\protect\citeauthoryear{Laplace}{Laplace}{1986}]{laplace1825essai}
Laplace, P. (1986).
\newblock {\em Essai philosophique sur les probabilit{\'e}s, r{\'e}{\'e}d}.
\newblock C. Bourgeois, Paris.

\bibitem[\protect\citeauthoryear{Lewis}{Lewis}{1973}]{lewis1973causation}
Lewis, D. (1973).
\newblock Causation.
\newblock {\em The Journal of Philosophy\/}~{\em 70\/}(17), 556--567.

\bibitem[\protect\citeauthoryear{Menzies}{Menzies}{2009}]{sep-causation-counte%
rfactual}
Menzies, P. (2009).
\newblock Counterfactual theories of causation.
\newblock {\em The Stanford Encyclopedia of Philosophy\/}.

\bibitem[\protect\citeauthoryear{Neyman}{Neyman}{1923}]{neyman1923application}
Neyman, J. (1923).
\newblock On the application of probability theory to agricultural experiments.
  essay on principles. section 9. translation of excerpts by d. dabrowska and
  t. speed.
\newblock {\em Statistical Science\/}~{\em 6\/}(1990), 462--47.

\bibitem[\protect\citeauthoryear{Nicholls, Tuzcu, Crowe, Sipahi, Schoenhagen,
  Kapadia, Hazen, Wun, Norton, Ntanios, and Nissen}{Nicholls
  et~al.}{2006}]{Nicholls2006}
Nicholls, S.~J., E.~M. Tuzcu, T.~Crowe, I.~Sipahi, P.~Schoenhagen, S.~Kapadia,
  S.~L. Hazen, C.-C. Wun, M.~Norton, F.~Ntanios, and S.~E. Nissen (2006, May).
\newblock {Relationship between cardiovascular risk factors and atherosclerotic
  disease burden measured by intravascular ultrasound.}
\newblock {\em Journal of the American College of Cardiology\/}~{\em 47\/}(10),
  1967--75.

\bibitem[\protect\citeauthoryear{Pearl}{Pearl}{2000}]{Pearl2000}
Pearl, J. (2000).
\newblock {\em {Causality : Models, reasoning, and inference}}.
\newblock Cambridge University Press.

\bibitem[\protect\citeauthoryear{Presanis, De~Angelis, Goubar, Gill, and
  Ades}{Presanis et~al.}{2011}]{presanis2011bayesian}
Presanis, A., D.~De~Angelis, A.~Goubar, O.~Gill, and A.~Ades (2011).
\newblock Bayesian evidence synthesis for a transmission dynamic model for hiv
  among men who have sex with men.
\newblock {\em Biostatistics\/}~{\em 12\/}(4), 666--681.

\bibitem[\protect\citeauthoryear{Ridker, Glynn, and Hennekens}{Ridker
  et~al.}{1998}]{Ridker1998}
Ridker, P., R.~Glynn, and C.~Hennekens (1998).
\newblock {C-reactive protein adds to the predictive value of total and HDL
  cholesterol in determining risk of first myocardial infarction}.
\newblock {\em Circulation\/}~{\em 97\/}(20), 2007.

\bibitem[\protect\citeauthoryear{R\o~ysland}{R\o~ysland}{2011}]{Røysland2011}
R\o~ysland, K. (2011).
\newblock {A martingale approach to continuous-time marginal structural
  models}.
\newblock {\em Bernoulli\/}~{\em 17\/}(3), 895--915.

\bibitem[\protect\citeauthoryear{Robins, Hern\'{a}n, and Brumback}{Robins
  et~al.}{2000}]{Robins2000}
Robins, J.~M., M.~a. Hern\'{a}n, and B.~Brumback (2000, September).
\newblock {Marginal structural models and causal inference in epidemiology.}
\newblock {\em Epidemiology (Cambridge, Mass.)\/}~{\em 11\/}(5), 550--60.

\bibitem[\protect\citeauthoryear{Rosenblum, Jewell, Laan, Shiboski, Straten,
  and Padian}{Rosenblum et~al.}{2009}]{rosenblum2009analysing}
Rosenblum, M., N.~Jewell, M.~Laan, S.~Shiboski, A.~Straten, and N.~Padian
  (2009).
\newblock Analysing direct effects in randomized trials with secondary
  interventions: an application to human immunodeficiency virus prevention
  trials.
\newblock {\em Journal of the Royal Statistical Society: Series A (Statistics
  in Society)\/}~{\em 172\/}(2), 443--465.

\bibitem[\protect\citeauthoryear{Rubin}{Rubin}{1978}]{Rubin1978}
Rubin, D. (1978).
\newblock {Bayesian inference for causal effects: The role of randomization}.
\newblock {\em The Annals of Statistics\/}, 34--58.

\bibitem[\protect\citeauthoryear{Sage and Melsa}{Sage and
  Melsa}{1971}]{sage1971system}
Sage, A. and J.~Melsa (1971).
\newblock {\em System identification}, Volume~80.
\newblock Academic Press.

\bibitem[\protect\citeauthoryear{{Van Gaal}, Mertens, and {De Block}}{{Van
  Gaal} et~al.}{2006}]{VanGaal2006}
{Van Gaal}, L.~F., I.~L. Mertens, and C.~E. {De Block} (2006, December).
\newblock {Mechanisms linking obesity with cardiovascular disease.}
\newblock {\em Nature\/}~{\em 444\/}(7121), 875--80.

\bibitem[\protect\citeauthoryear{Whitmore}{Whitmore}{1986}]{Whitmore1986}
Whitmore, G.~A. (1986).
\newblock {Normal-gamma mixtures of inverse Gaussian distributions}.
\newblock {\em Scandinavian Journal of Statistics\/}~{\em 13\/}(3), 211--220.

\bibitem[\protect\citeauthoryear{Wimsatt}{Wimsatt}{1994}]{Wimsatt1994}
Wimsatt, W. (1994).
\newblock {The ontology of complex systems: levels of organization,
  perspectives, and causal thickets}.
\newblock {\em Canadian Journal of Philosophy\/}~{\em 20}, 207--274.

\bibitem[\protect\citeauthoryear{Wimsatt}{Wimsatt}{2007}]{Wimsatt2007}
Wimsatt, W.~C. (2007, February).
\newblock {Aggregate, composed, and evolved systems: Reductionistic heuristics
  as means to more holistic theories}.
\newblock {\em Biology \& Philosophy\/}~{\em 21\/}(5), 667--702.

\bibitem[\protect\citeauthoryear{Zhu, Song, and Taylor}{Zhu
  et~al.}{2011}]{zhu2011stochastic}
Zhu, B., P.~Song, and J.~Taylor (2011).
\newblock Stochastic functional data analysis: A diffusion model-based
  approach.
\newblock {\em Biometrics\/}.

\end{thebibliography}
\end{document}